\documentclass[10pt,twoside]{article}
\usepackage{Latex-document}

\newtheorem{thm}{Theorem}
\newtheorem{prop}{Proposition}

\title{\bf Vector Bundles on a K3 Surface\thanks{Supported in part by the JSPS Grant-in-Aid for Scientific Research
(A) (2) 10304001.}\vskip 6mm}
\author{Shigeru Mukai\vspace*{-0.5cm}\thanks{Research Institute for Mathematical Sciences,
Kyoto University, Kyoto 606-8502, Japan. E-mail: mukai@kurims.kyoto-u.ac.jp}}
\date{\vspace{-8mm}}

\begin{document}

\maketitle

\thispagestyle{first} \setcounter{page}{495}

\begin{abstract}

\vskip 3mm

A K3 surface is a quaternionic analogue of an elliptic curve from a view point of moduli of vector bundles. We can
prove the algebraicity of certain Hodge cycles and a rigidity of curve of genus eleven and gives two kind of
descriptions of Fano threefolds as applications. In the final section we discuss a simplified construction of
moduli spaces.

\vskip 4.5mm


\end{abstract}

\vskip 12mm

\section{Introduction}

\vskip-5mm \hspace{5mm}

A locally free sheaf $E$ of ${\cal O}_X$-modules is called {\it a
vector bundle} on  an algebraic variety  $X$. As a natural
generalization of line bundles vector bundles have two important
roles in algebraic geometry. One is the linear system. If  $E$  is
generated by its global sections $H^0(X, E)$, then it gives rise
to a morphism to a Grassmann variety, which we denote by $\Phi_E :
X \longrightarrow G(H^0(E), r)$, where $r$  is the rank of  $E$.
This morphism is related with the classical linear system by the
following diagram:
\begin{equation}
\begin{array}{rcccl}
  & X           & \stackrel{\Phi_E}\longrightarrow  & G(H^0(E), r) & \\
\Phi_L &\downarrow   &         & \downarrow & {\rm Pl\ddot{u}cker} \\
&{\bf P}^*H^0(L) & \cdots\rightarrow  & {\bf P}^*(\bigwedge^r
H^0(E)), &
\end{array}
\end{equation}
where  $L$  is the determinant line bundle of $E$ and $\Phi_L$ is
the morphism associated to it.

The other role is the moduli. The moduli space of line bundles
relates a (smooth complete) algebraic curve with an abelian
variety called the {\it Jacobian variety}, which is crucial in the
classical theory of algebraic functions in one variable. The
moduli of vector bundles also gives connections among different
types varieties, and often yields new varieties that are difficult
to describe by other means.

In higher rank case it is natural to consider the moduli problem
of  $E$ under the restriction that $\det E$  is unchanged. In view
of the above diagram, vector bundles and their moduli reflect the
geometry of the morphism  $X \longrightarrow {\bf P}^*H^0(L)$ via
Grassmannians and Pl\"ucker relations. In this article we consider
the case where  $X$ is a K3 surface, which is one of two
2-dimensional analogues of an elliptic curve and  seems an ideal
place to see such reflection.

\section{Curves of genus one}

\vskip-5mm \hspace{5mm}

The moduli space of line bundles on an algebraic variety is called
the {\it Picard variety}. The Picard variety  ${\rm Pic}\, C$ of
an algebraic curve  $C$ is decomposed into the disjoint union
$\coprod_{d \in {\bf Z}}{\rm Pic}_d\, C$ by the degree $d$ of line
bundles. Here we consider the case of genus 1. All components
${\rm Pic}_d\, C$ are isomorphic to $C$ if the ground field is
algebraically closed. \footnote{More precisely, this holds true if
$C$  has a rational point.} But this is no more true otherwise.
For example the Jacobian ${\rm Pic}_0\, C$ has always a rational
point but $C$ itself does not. \footnote{Two components  ${\rm
Pic}_0\, C$ and  ${\rm Pic}_{g-1}\, C$ deserve the name {\it
Jacobian}. They coincide in our case $g=1$.} We give other
examples:

\smallskip\noindent
{\bf Example 1}\, Let  $C_4$  be an intersection of two quadrics
$q_1(x) = q_2(x) = 0$ in the projective space ${\bf P}^3$ and $P$
the pencil of defining quadrics. Then the Picard variety ${\rm
Pic}_2\, C_4$ is the double cover  of $P \simeq {\bf P}^1$ and the
branch locus consists of 4 singular quadrics in $P$. Precisely
speaking, its equation is given by  $\tau^2 = {\rm disc}\,
(\lambda_1q_1 + \lambda_2q_2)$.

\smallskip

Let  $G(2, 5) \subset {\bf P}^9$ be the 6-dimensional Grassmann
variety embedded into ${\bf P}^9$ by the Pl\"ucker coordinate. Its
projective dual is the dual Grassmannian $G(5, 2) \subset \hat
{\bf P}^9$, where $G(2,5)$ parameterizes 2-dimensional subspaces
and $G(5,2)$ quotient spaces.

\smallskip\noindent
{\bf Example 2}\, A transversal linear section  $C = G(2, 5) \cap
H_1 \cap \cdots \cap H_5$ is a curve genus 1 and of degree 5. Its
Picard variety  ${\rm Pic}_2\, C$ is isomorphic to the dual linear
section $\hat C = G(5, 2) \cap \langle H_1, \ldots, H_5\rangle$,
the intersection with the linear subspace spanned by 5 points $
H_1, \ldots, H_5 \in \hat {\bf P}^9$.
\smallskip

\section{Moduli K3 surfaces}

\vskip-5mm \hspace{5mm}

A compact complex 2-dimensional manifold $S$  is {\it a K3
surface} if the canonical bundle is trivial and the irregularity
vanishes, that is, $K_S= H^1({\cal O}_S)=0$. A smooth quartic
surface  $S_4 \subset {\bf P}^3$ is the most familiar example. Let
us first look at the 2-dimensional generalization of Example 1:

\smallskip\noindent
{\bf Example 3}\, Let  $S_8$  be an intersection of three general
quadrics in ${\bf P}^5$ and
 $N$ the net of defining quadrics.
Then the moduli space $M_S(2, {\cal O}_S(1), 2)$ is a double cover
of  $N \simeq {\bf P}^2$ and the branch locus, which is of degree
6, consists of singular quadrics in $N$.

\smallskip
Here $M_S(r, L, s)$, $L$ being a line bundle, is the moduli space
of stable sheaves  $E$  on a K3 surface $S$ with rank $r$, $\det E
\simeq L$ and  $\chi(E) = r+s$. Surprisingly two surfaces $S_8 =
(2) \cap (2) \cap (2) \subset {\bf P}^5$ and  $M_S(2,{\cal
O}_S(1), 2) \stackrel{2:1}\longrightarrow {\bf P}^2$ in this
example are both K3 surfaces. This is not an accident. In respect
of moduli space, vector bundles a K3 surface look like Picard
varieties in the preceding section.

\begin{thm}{\rm (\cite{Mu},\cite{tata})}
The moduli space $M_S(r, L, s)$ is smooth of dimension $(L^2) -2rs
+2$. $M_S(r, L, s)$ is again a K3 surface if it is compact and and
of dimension $2$.
\end{thm}

A K3 surface $S$ and a moduli K3 surface appearing as $M_S(r, L,
s)$ are not isomorphic in general \footnote{We take the complex
number field {\bf C} as ground field except for sections 2 and 7.}
but their polarized Hodge structures, or periods, are isomorphic
to each other over ${\bf Q}$ (\cite{tata}). The moduli is not
always fine but there always exists a universal ${\bf
P}^{r-1}$-bundle over the product $S \times M_S(r, L, s)$. Let
${\cal A}$  be the associated sheaf of Azumaya algebras, which is
of rank $r^2$ and locally isomorphic to the matrix algebra
$Mat_r({\cal O}_{S \times M})$. ${\cal A}$ is isomorphic to ${\cal
E}nd\, {\cal E}$ if a universal family ${\cal E}$ exists. The
Hodge isometry between  $H^2(S, {\bf Q})$ and $H^2(M_S(r, L, s),
{\bf Q})$ is given by
$c_2({\cal A})/2r \in H^4(S \times M, {\bf Q}) \simeq H^2(S, {\bf
Q})^\vee \otimes H^2(M, {\bf Q})$.

Example 2 has also a K3 analogue. Let $\Sigma_{12} \subset {\bf
P}^{15}$ be a 10-dimensional spinor variety $SO(10)/U(5)$, that
is, the orbit of a highest weight vector in the projectivization
of the 16-dimensional spinor representation. The anti-canonical
class is 8 times the hyperplane section and a transversal linear
section  $S = \Sigma_{12} \cap H_1 \cap \cdots \cap H_8$ is a K3
surface (of degree 12). As is similar to $G(2,5)$ the projective
dual $\hat\Sigma_{12} \subset \hat{\bf P}^{15}$ of $\Sigma_{12}$
is again a 10-dimensional spinor variety.

\smallskip\noindent
{\bf Example 4}\, The moduli space $M_S(2, {\cal O}_S(1), 3)$  is
isomorphic to the dual linear section $\hat S = \hat\Sigma_{12}
\cap \langle H_1, \ldots, H_8\rangle$.
\smallskip

In this case, moduli is fine and the relation between $S$ and the
moduli K3 are deeper. The universal vector bundle ${\cal E}$ on
the product gives an equivalence between the derived categories
${\bf D}(S)$ and ${\bf D}(\hat S)$ of coherent sheaves, the
duality $\hat{\hat S} \simeq S$ holds (cf. \cite{dual}) and
moreover the Hilbert schemes ${\rm Hilb}^2\, S$ and ${\rm
Hilb}^2\, \hat S$ are isomorphic to each other.

\smallskip\noindent
{\bf Remark}\, (1) Theorem 1 is generalized to the non-compact
case by Abe \cite{Ab}.

\noindent (2) If a universal family ${\cal E}$ exists, the derived
functor with kernel ${\cal E}$ gives an equivalence of derived
categories of coherent sheaves on $S$ and the moduli K3
(Bridgeland \cite{B}). In even non-fine case the derived category
${\bf D}(S)$ is equivalent to that of the moduli K3 $M$ twisted by
a certain element $\alpha \in H^2(M, {\bf Z}/r{\bf Z})$ (C\v ald\v
araru \cite{C}).

\section{Shafarevich conjecture}

\vskip-5mm \hspace{5mm}

Let  $S$ and  $T$ be algebraic K3 surfaces and $f$ a Hodge
isometry between $H^2(S, {\bf Z})$ and $H^2(T, {\bf Z})$. Then the
associated cycle $Z_f \in H^4(S \times T, {\bf Z}) \simeq H^2(S,
{\bf Z})^\vee \otimes H^2(T, {\bf Q})$ on the product $S \times T$
is a Hodge cycle. This is algebraic by virtue of the Torelli type
theorem of Shafarevich and Piatetskij-Shapiro. Shafarevich
conjectured in \cite{Sh} a generalization to Hodge isometries over
${\bf Q}$. Our moduli theory is able to answer this affirmatively.

\begin{thm}
Let $f: H^2(S, {\bf Q}) \longrightarrow H^2(T, {\bf Q})$ be a
Hodge isometry. Then the associated (Hodge) cycle $Z_f \in H^4(S
\times T, {\bf Q})$ is algebraic.
\end{thm}

In \cite{tata}, we already proved this partially using Theorem 1
(cf. \cite{Ni} also). What we need further is the moduli space of
projective bundles. Let $P \longrightarrow S$ be a ${\bf
P}^{r-1}$-bundle over  $S$. The cohomology class $[P] \in H^1(S,
PGL(r, {\cal O}_S))$ and the natural exact sequence (in the
classical topology)
$$0 \longrightarrow {\bf Z}/r{\bf Z} \longrightarrow SL(r, {\cal O}_S)
\longrightarrow PGL(r, {\cal O}_S) \longrightarrow 1$$ define an
element of $H^2(S, {\bf Z}/r{\bf Z})$, which we denote by $w(P)$.

Fix $\alpha \in H^2(S, {\bf Z})$ and $r$, we consider the moduli
of ${\bf P}^{r-1}$-bundles $P$ over  $S$ with $w(P) = \alpha$ mod
$r$ which are stable in a certain sense. If the self intersection
number  $(\alpha^2)$ is divisible by $2r$, then the moduli space
contains a 2-dimensional component, which we denote by
$M_S(\alpha/r)$. The following, a honest generalization of
computations in \cite{tata}, is the key of our proof:

\begin{prop}
Assume that $(\alpha^2)$ is divisible by $2r^2$. Then
$H^2(M_S(\alpha/r), {\bf Z})$  is isomorphic to $L_0 + {\bf
Z}\alpha/r \subset  H^2(S, {\bf Q})$ as polarized Hodge structure,
where  $L_0$  is the submodule of $H^2(S, {\bf Z})$ consisting of
$\beta$ such that the intersection number $(\beta. \alpha)$  is
divisible by $r$.
\end{prop}

For example  let $S_2$ be a double cover of ${\bf P}^2$ with
branch sextic. If $\alpha \in H^2(S, {\bf Z})$ satisfies $(\alpha.
h) \equiv 1  \,{\rm mod}\, 2$ and $(\alpha^2) \equiv 0  \,{\rm
mod}\, 4$, then $M_S(\alpha/2)$  is a K3 surface of degree 8. This
is the inverse correspondence of Example 1 (cf. \cite{Tj},
\cite{Ne}). Details will be published elsewhere.

\section{Non-Abelian Brill-Noether locus}

\vskip-5mm \hspace{5mm}

Let $C$ be a smooth complete algebraic curve. As a set a
Brill-Noether locus of  $C$  is a stratum of the Picard variety
${\rm Pic}\, C$ defined by $h^0(L)$,  the number of global
sections of a line bundle  $L$. The standard notation is

$$W^r_d = \{ [L] \,|\, h^0(L) \ge r + 1\} \subset {\rm Pic}_d\, C,$$
for which we refer \cite{ACGH}. Non-Abelian analogues are defined
in the moduli space  ${\cal U}_C(2)$ of stable 2-bundles on $C$
similarly. The {\it non-Abelian Brill-Noether locus of type III}
is

$${\cal SU}_C(2, K:n) = \{F \,|\, \det F \simeq {\cal O}_C(K_C), \, h^0(F)
\ge n \} \subset {\cal U}_C(2)$$ for a non-positive integer $n$,
and {\it type II} is

$${\cal SU}_C(2, K:nG) = \{F \,|\, \det F \simeq \det G \otimes {\cal
O}_C(K_C), \, \dim {\rm Hom}(G, F) \ge n \} \subset {\cal
U}_C(2)$$ for a vector bundle  $G$  of rank 2 and $n \equiv \deg G
\,{\rm mod}\, 2$. By virtue of the (Serre) duality, these have
very special determinantal descriptions. We give them scheme
structures using these descriptions (\cite{BN}).

Assume that  $C$ lies on a K3 surface  $S$. If $E$  belongs to
$M_S(r, L, s)$, then the restriction  $E|_C$ is of canonical
determinant and we have $h^0(E|_C) \ge h^0(E) \ge \chi(E)= r+s$.
So $E|_C$ belongs to ${\cal SU}_C(2, K:r+s)$ if it is stable. This
is one motivation of the above definition. The case of genus 11,
the gap value of genera where Fano 3-folds of the next section do
not exist, is the molst interesting.

\begin{thm}{\rm (\cite{eleven})}
If  $C$  is a general curve of genus $11$, then the Brill-Noether
locus $T = {\cal SU}_C(2, K:7)$ of type III is a K3 surface and
the restriction $L$ of the determinant line bundle is of degree
$20$.
\end{thm}

There exists a universal family  ${\cal E}$ on $C \times T$. We
moreover have the following:

\begin{itemize}
\item the restriction ${\cal E}|_{x \times T}$ is is stable and belongs to
$M_T(2,L, 5)$, for every  $x \in C$, and

\item the classification morphism $C \longrightarrow \hat T = M_T(2, L, 5)$
is an embedding.
\end{itemize}
This embedding is a non-Abelian analogue of the Albanese morphism $X \longrightarrow$ \linebreak ${\rm Pic}_0
({\rm Pic}_0 X)$ and we have the following:

\noindent {\bf Corollary}\it \, A general curve of genus eleven
has a unique embedding to a K3 surface. \rm

In \cite{MM}, we studied the forgetful map $\varphi_g$  from the
moduli space ${\cal P}_g$ of pairs  of a curve $C$  of genus  $g$
and a K3 surface $S$ with $C \subset S$ to the moduli space ${\cal
M}_g$ of curves of genus  $g$ and the generically finiteness of
$\varphi_{11}$. The above correspondence  $C \mapsto \hat T$ gives
the inverse rational map of $\varphi_{11}$. We recall the fact
that $\varphi_{10}$ is not dominant in spite of the inequality
$\dim {\cal P}_{10} = 29 > \dim {\cal M}_{10}=27$ (\cite{nagata}).

\section{Fano 3-folds}

\vskip-5mm \hspace{5mm}

A smooth 3-dimensional projective variety is called a {\it Fano
3-fold} if the anti-canonical class $-K_X$  is ample. In this
section we assume that the Picard group  ${\rm Pic}\, X$  is
generated by $-K_X$. The self intersection number $(-K_X)^3 =
2g-2$  is always even and the integer  $g \ge 2$  is called the
{\it genus}, by which the Fano 3-folds are classified into 10
deformation types. The values of $g$ is equal to $2, \ldots, 10$
and $12$. A Fano 3-folds of genus $g \le 5$  is a complete
intersection of hypersurfaces in a suitable weighted projective
space.

By Shokurov \cite{Sho}, the anticanonical linear system $|-K_X|$
always contain a smooth member $S$, which is a K3 surface. In
\cite{PNAS} we classified Fano 3-folds $X$ of Picard number one
using {\it rigid bundles}, that is, $E \in M_S(r, L, s)$ with
$(L^2) -2rs = -2$. For example  $X$ is isomorphic to a linear
section of the 10-dimensional spinor variety, that is,
\begin{equation}\label{linsec}
X \simeq \Sigma_{12} \cap H_1 \cap \cdots \cap H_7,
\end{equation}
in the case of genus 7 and a linear section
\begin{equation}\label{linsec2}
X \simeq \Sigma_{16} \cap H_1 \cap H_2 \cap H_3,
\end{equation}
of the 6-dimensional symplectic, or Lagrangian, Grassmann variety $\Sigma_{16}=SP(6)/$ $U(3) \subset {\bf P}^{13}$
in the case of genus 9. The non-Abelian Brill-Noether loci shed new light on this classification.
\begin{thm}
A Fano 3-fold  $X$  of genus $7$ is isomorphic to  the
Brill--Noether locus ${\cal SU}_C(2, K:5)$  of Type~III for a
smooth curve  $C$  of genus $7$.
\end{thm}

This description is dual to the description (\ref{linsec}) in the
following sense. First two ambient spaces of  $X$, the moduli
${\cal U}_C(2)$ and the Grassmannian $G(5, 10) \supset
\Sigma_{12}$ are of the same dimension. Secondly let $N_1$  and
$N_2$  be the normal bundles of  $X$ in these ambient spaces. Then
we have $N_1 \simeq N_2^\vee \otimes {\cal O}_X(-K_X)$, that is,
two normal bundles are twisted dual to each other.

 \begin{thm}
A Fano $3$-fold of genus $9$ is isomorphic to the Brill--Noether
locus  ${\cal SU}_C(2,K\colon 3G)$ of Type~II for a nonsingular
plane quartic curve  $C$ and $G$ a rank $2$ stable vector bundle
over $C$ of odd degree.
 \end{thm}

This descriptions is also dual to (\ref{linsec2}) in the above
sense: The moduli  ${\cal U}_C(2)$ and the Grassmannian $G(3, 6)
\supset \Sigma_{16}$ are of the same dimension and the two normal
bundles of $X$ are twisted dual to each other. Each Fano 3-fold of
genus 8, 12 and conjecturally 10 has also such a pair of
descriptions.

\section{Elementary construction}

\vskip-5mm \hspace{5mm}

The four examples in sections 1 and 2  are very simple and invite
us to a simplification of moduli construction. Let $C_4$ be as in
Example 1 and ${\cal M}at_2$ the affine space associated to the
16-dimensional vector space $\bigoplus_{i=0}^3 ({\bf C}^2
\otimes{\bf C}^2)x_i$, where $(x_i)$ is the homogeneous coordinate
of ${\bf P}^3$. Let ${\cal M}at_{2,1}$ be the closed subscheme
defined by the condition that

\smallskip
\hskip 0.9cm$A(x) = \sum_{i=0}^3 A_ix_i \in {\cal M}at_2$  is of
rank $\le 1$ everywhere on $C_4$

\smallskip\noindent
and $R$ its coordinat ring. Then the Picard variety ${\rm Pic}_2\,
C_4$  is the projective spectrum ${\rm Proj}\, R^{SL(2)\times
SL(2)}$ of the invariant ring by construction. (See \cite{moduli}
for details.) The above condition is equivalent to that $\det
A(x)$  is a linear combination of $q_1(x)$ and  $q_2(x)$. The
invariant ring is generated by three elements by Theorem 2.9.A  of
Weyl \cite{Wey}. Two of them, say $B_1, B_2$, are of degree 2 and
correspond to $q_1(x)$ and  $q_2(x)$, respectively. The rest, say
$T$ of degree 4, is the determinant of 4 by 4 matrix obtained from
the four coefficients $A_0, \ldots, A_3 \in {\bf C}^2 \otimes {\bf
C}^2$ of $A(x)$. There is one relation $T^2 = f_4(B_1, B_2)$ .
Hence ${\rm Proj}\, R^{SL(2)\times SL(2)}$ is a double cover of
${\bf P}^1$ as desired.

The moduli space  $M_S(2, {\cal O}_S(1), 2)$ in Example 3 is
constructed similarly. Let ${\cal A}lt_4$ be the affine space
associated to the vector space $\bigoplus_{i=0}^5 (\bigwedge^2{\bf
C}^4)x_i$ and ${\cal A}lt_{4,2}$ the subscheme defined by the
condition that $\sum_{i=0}^5 A_ix_i \in {\cal A}lt_4$ is of rank
$\le 2$ everywhere on $S_8$. Then the invariant ring of the action
of $SL(4)$ on  ${\cal A}lt_{4,2}$ is generated by four elements
$B_1, B_2, B_3, T$  of degree 2, 2, 2, 6. There is one relation
$T^2 = f_6(B_1, B_2, B_3)$ and $M_S(2,{\cal O}_S(1), 2)$, the
projective spectrum ${\rm Proj}\, R^{SL(4)}$, is a double cover of
${\bf P}^2$ as described.

The moduli space of vector bundles on a surface was first
constructed by Gieseker \cite{Gi}. He took the Mumford's GIT
quotient \cite{GIT} of Grothendieck's Quot scheme \cite{Gr} by
$PGL$ and used the Gieseker matrix to measure the stability of the
action. In the above construction, we take the quotient of ${\cal
A}lt_{4,2}$, which is nothing but the affine variety of Gieseker
matrices of suitable rank 2 vector bundles, by a general linear
group $GL(4)$.

The Jacobian, or the Picard variety, of a curve is more
fundamental. Weil \cite{Wei} constructed ${\rm Pic}_g\, C$ as an
algebraic variety using the symmetric product ${\rm Sym}^g\, C$
and showed its projectivity by Lefschetz' $3\Theta$ theorem. Later
Seshadri and Oda \cite{SO} constructed ${\rm Pic}_d\, C$ for
arbitray $d$ (over the same ground field as $C$) by also taking
the GIT quotient of Quot schemes. The above constructions
eliminate  Quot schemes and the concept of linearization from
those of Gieseker, Seshadri and Oda.

\label{lastpage}

\end{document}